\documentclass[11pt]{amsart}
\usepackage[dvips]{graphicx}
\usepackage{amsmath, amsthm, amscd, amssymb, color}
\usepackage{epsfig}
\usepackage[all]{xy}

\newtheorem{thm}{Theorem}[section]
\newtheorem*{thm*}{Theorem}
\newtheorem{cor}[thm]{Corollary} 
\newtheorem*{cor*}{Corollary}
\newtheorem{lem}[thm]{Lemma}

\newtheorem*{con*}{Conjecture}
\newtheorem*{prob*}{Problem}

\theoremstyle{definition}
\newtheorem{defn}[thm]{Definition}

\theoremstyle{remark}
\newtheorem{rem}[thm]{Remark}

\newcommand{\N}{\mathcal{N}}
\newcommand{\F}{\mathcal{F}}
\newcommand{\p}{\mathcal{P}}

\newcommand{\s}{\mathcal{S}}

\newcommand{\bbC}{\mathbb{C}}

\newcommand{\bbR}{\mathbb{R}}

\newcommand{\bbN}{\mathbb{N}}

\begin{document}
\title{Zeta Functions of Infinite Graph Bundles}

\author[Samuel Cooper]{Samuel Cooper$^{*}$}
\address{Department of Mathematics,
Vanderbilt University, Nashville, TN 37240}
\email{samuel.d.cooper@vanderbilt.edu}
\thanks{$^{*}$Partially supported by an NSF REU grant}

\author[Stratos Prassidis]{Stratos Prassidis$^{**}$}
\address{Department of Mathematics
Canisius College, Buffalo, NY 14208}
\email{prasside@canisius.edu}
\thanks{$^{**}$Partially supported by Canisius College Summer
Grant and an NSF REU grant}

\begin{abstract}
We compute the equivariant zeta function for bundles over infinite
graphs and for infinite covers. 
In particular, we give a ``transfer formula'' for the zeta function of 
infinite graph covers. Also, when the infinite cover is given as a limit
of finite covers, we give a formula for the limit of the zeta
functions.
\end{abstract}

\subjclass[2000]{Primary 11M36; Secondary 05C50, 11M41}

\maketitle

\section{Introduction}

The Ihara zeta function of a finite graph reflects combinatorial and
spectral properties of that graph (\cite{ihara}, \cite{bass}, \cite{stark1}). 
Originally, Ihara defined the zeta function on finite graphs imitating
the classical definition of the zeta function: 
$${\zeta}_X(z) = \prod_{[C]} (1-z^{\ell(C)})^{-1},$$ 
where the product is over all equivalence classes of primitive closed loops 
$C$ in $X$ and $\ell(C)$ denotes the length of $C$.  
In \cite{bass}, it was 
shown that, for a finite graph $X$:
$${\zeta}_X(z)^{-1} = (1- z^2)^{\epsilon - \nu} det(I - zA + z^2Q),$$
where $\epsilon$ is the number of edges, $\nu$ is the number of vertices,
$A$ is the adjacency matrix of $X$, and $Q$ is the diagonal matrix with
entries $\text{deg}(v) - 1$, for each $v\in V(X)$.
In \cite{grig}, the definition of Ihara zeta function was extended to
infinite graphs that are limits of sequences of finite graphs. In particular,
it was shown in \cite{grig}, using the results in \cite{serre}, 
that the sequence of the zeta functions of the finite graphs converges.
In \cite{clair}, \cite{clair2}, \cite{guido1}, \cite{guido2}, \cite{guido3},
the expression of the zeta function as a rational function was extended to
infinite graphs that admit an action of a discrete group $\Gamma$ with
finite quotient. The determinant in the finite case is replaced by
the determinant in a von Neumann algebra ${\N}_0(X, {\Gamma})$ of 
all the bounded operators on $L^2(V(X))$. 
In \cite{feng} the zeta function of finite graph bundles over 
finite was computed generalizing the results on 
graph coverings that appear \cite{stark1}, \cite{stark2},
\cite{stark3}. Their results can be described as transfer 
results for the Ihara zeta function.

We combine the results on infinite graphs and bundles to derive
a transfer formula for infinite bundles and coverings.
Let $\phi$ be an $\text{Aut}(F)$-assignment on $X$. Let $({\Gamma},
{\Delta})$ be a pair of groups that act on $X$ and $F$ in such a way
that the actions are $\phi$-compatible and by finite co-volume.

\begin{thm*}[Main Theorem 1]
With the above assumptions, the equivariant zeta function is given by:
$${\zeta}_{X \times^{\psi}F,{\Gamma}{\times}{\Delta}}(z)^{-1} = 
(1 - z^2)^{-\chi^{(2)}(X \times^{\psi} F)} \; 
{\text{det}}_{{\Gamma}{\times}{\Delta}} 
\left(I - \sum_{\gamma \in 
\text{Aut}(F)}(A_{{\overrightarrow{X}}_{(\psi, \gamma)}} \otimes P_\gamma + 
I_X \otimes A_F)z + Qz^2\right),$$
where $A_{{\overrightarrow{X}}_{(\phi, \gamma)}}$ is the 
adjacency matrix of the directed graph spanned by the edges in 
${\psi}^{-1}({\gamma})$, $P_{\gamma}$ is the permutation matrix induced by 
the action of $\gamma$  on $V(F)$, ${\chi}^{(2)}$ is the Euler characteristic 
of the quotient $X_1{\times}^{{\phi}_1}F$, $\text{det}_{{\Gamma}{\times}
{\Delta}}$ is 
the determinant defined on the von Neumann algebra of ${\Gamma}{\times}
{\Delta}$, and $Q$ is
the diagonal operator such that $Q(y, i) = \text{deg}(y) + \text{deg}(i) - 1$.
\end{thm*}

Using similar methods, we prove a decomposition formula of the Ihara zeta 
function for infinite covers. Let $p: Y \to X$ be a cover with $X$ finite.
Let ${\Gamma} = \text{Cov}(p)$.

\begin{thm*}[Main Theorem 2] 
With the above notation,
$${\zeta}_{Y,{\Gamma}}(z)^{-1} = (1 - z^2)^{-\chi^{(2)}(Y)} \; 
{\text{det}}_{\Gamma} \left(I - \left(\sum_{\gamma \in \Gamma} 
A_{{({\overrightarrow{Y}}_{(\psi, \gamma)}})} \otimes 
P_\gamma\right)z + Qz^2 \right),$$ 
where $Q$ is the diagonal operator with $(x, {\gamma})$-entry
$\text{deg}(x) - 1$. 
\end{thm*}

We apply the above calculations to sequences of {\it
strongly convergent graphs}. In particular, a sequence 
$\{(X_n, w_n)\}_{n\in\bbN}$ is strongly convergent to $(X, w)$ if it is a
covering sequence of regular graphs converging to $X$ in such a way 
that $X$ covers compatibly each element of the sequence. Such sequences
appear when we consider the Cayley graphs of finite quotients of a groups
converging to the Cayley graph of the group.

\section{Preliminaries}

We now define a number of terms that we will use later on.

\begin{defn} 
Let $G$ be any locally finite graph. 
Then we define the adjacency operator 
$A_G$ as follows: for any $u, v \in V(G)$,
$$A_G(u, v) = 
\begin{cases}
1, & \text{ if } u \sim v \\
0, & \text{ otherwise. } 
\end{cases}$$
\end{defn}

The definition makes sense even if the graph is directed. If $G$ is
undirected, the $A_G$ is symmetric.

\begin{defn}
Let $\tilde{G}$ and $G$ be locally finite graphs.  We say that 
$$p: \tilde{G} \to G$$ 
is a graph covering if the following two conditions hold:
\begin{enumerate}
\item If $x \sim_{\tilde{G}} y$, then $p(x) \sim_G p(y)$.
\item For any $x \in \tilde{G}$, $p: N(x) \to N\big(p(x)\big)$ is a bijection.
\end{enumerate}
\end{defn} 

The first condition in the definition means that $p$ is a graph map. The
second condition is a local triviality condition.

Graph bundles are defined in \cite{kwak}. They generalize the graph coverings
in the sense that the ``fiber graph'' is allowed to have a non-empty set
of edges. We will concentrate on bundles with finite fibers. For a graph
$X$, we denote by $E(\overrightarrow{X})$ the set of ordered edges--i.e.,
each edge of $X$ appears twice, each with opposite orientation.

\begin{defn}
Let $G$ be any locally finite graph (possibly infinite), let $F$ be a 
finite graph.  We define an $\text{Aut}(F)$-voltage assignment on $G$ by  
$$\phi: E({\overrightarrow{G}}) \to \text{Aut}(F), 
\quad \phi(uv) = \phi(vu)^{-1}.$$ 
\end{defn}

\begin{defn} Let $G$ be a locally finite graph, $F$ a finite graph, 
and $\phi$ an Aut($F$)-voltage assignment on $G$.  We define a graph 
bundle $G \times^\phi F$ to be the graph with vertex set 
$V(G) \times V(F)$, with two vertices $(u, i), (v, j) \in G \times^\phi F$ 
adjacent if either one of the following two conditions hold: 
\begin{enumerate}
\item $ u \sim v \text{ and } j = i^{\phi(uv)} $
\item $ u = v \text{ and } i \sim j. $
\end{enumerate}
\end{defn}

Let $\phi$ be a $\text{Aut}(F)$-voltage assignment on $G$. Let ${\gamma} \in
 \text{Aut}(F)$. 
\begin{enumerate}
\item Let 
${{\overrightarrow{G}}}_{(\phi, \gamma)}$ denote the spanning subgraph of 
the digraph ${\overrightarrow{G}}$ whose directed edge set is 
$\phi^{-1}(\gamma)$.
\item We define the permutation operator $P_{\gamma}$ by the following 
formula: for any two vertices $i, j$  in $V(F)$, 
$$P_\gamma (i, j) = 
\begin{cases}
1, & \text{ if } j = i^\gamma \\
0, & \text{ otherwise. }
\end{cases}$$ 
\end{enumerate}

\begin{rem}
When the graphs are infinite, the matrices defined above are operators
on the Hilbert space with basis the vertex set of the graph. More
precisely, if $G$ is any locally finite graph, set $L^2(G)$ to be
the Hilbert space:
$$L^2(G) = \left\{f: V(G) \to \bbC :\;\; \sum_{u\in V(G)}|f(u)|^2 < 
\infty \right\}.$$
Then the adjacency operator is given by
$$A(f)(u) = \sum_{u \sim v}f(v).$$ 
With the above notation,
$$P_{\gamma}(f)(i) = f({\gamma}(i)).$$
\end{rem}

The following combines covering maps and bundles. 

\begin{thm}{\label{cover}}
Let $F$ and $X$ be locally finite graphs. Let $X$ be equipped with
an $\text{Aut}(F)$-voltage assignment $\phi$. Let $p: Y \to X$
be a covering map and $\psi$ the $\text{Aut}(F)$-voltage assignment
$$\psi : E(\overrightarrow{Y}) \to \text{Aut}(F), \quad
{\psi}(xy) = {\phi}(p(x)p(y)).$$
Define a graph map
$$\tilde{p}: Y \times^{\psi} F \to X{\times}^{\phi}F, \quad
\tilde{p}(x, i) = (p(x), i).$$
Then $\tilde{p}$ is a covering map.
\end{thm}

\begin{proof}
First we will show prove that $\tilde{p}$ is a graph map i.e., that is
preserves adjacency.  Let $(x, i) \sim (y, j)$ in 
$Y{\times}^{\psi}F$.  There are two cases to consider: 
\begin{enumerate}
\item Suppose $x \sim y$ in $Y$. Then $p(x) \sim p(y)$, and  
$j = i^{{\psi}(xy)} = i^{\phi(p(x)p(y))}$. Thus, by definition, 
$$\tilde{p}(x, i) = (p(x), i) \sim (p(y), j) =\tilde{p}(y, j), \;\;
\text{in}\;\; X{\times}^{\phi}F.$$
\item Suppose $x = y$. Then $i \sim j$ in $F$, and clearly $p(x) = p(y)$.
Thus, by definition, $\tilde{p}(x, i) \sim \tilde{p}(y, j)$.
\end{enumerate}
Thus $\tilde{p}$ preserves adjacency.

Now we must show that $\tilde{p}|_{N(x, i)}$ is a bijection.  

\vspace{12pt}
\noindent
\underline{ $\tilde{p}|_{N(x, i)}$ is an injection}.
Let $(y_1, j_1), (y_2, j_2) \in N(x, i)$ with 
$\tilde{p}(y_1, j_1) = \tilde{p}(y_2, j_2).$  Then we know that 
$p(y_1) = p(y_2)$ and $j_1 = j_2.$  Now there are two cases to consider:
\begin{enumerate}
\item Suppose $y_1 = x$ and $i \sim j_1$. Then $i \sim j_2$, and since 
$(x, i) \sim (y_2, j_2)$, we must have $x = y_2$. Thus $y_1 = y_2$, so 
$(y_1, j_1) = (y_2, j_2)$. The same argument works if $y_2 = x$.
\item  Suppose $y_1 \sim x$, and $i = {j_1}^{{\psi}(y_1 x)} = 
{j_1}^{\phi(p(y_1)p(x))}$. Since $p(y_1) = p(y_2)$, we see that 
$i = {j_2}^{\phi(p(y_2)p(x))}$. Thus, since $(y_2, j_2) \sim (x, i)$, we must 
have $y_2 \sim x$. Now, since $p$ is a graph covering map,  
$p|_{N(x)}$ is a bijection. But $y_1, y_2 \in N(x)$ and $p(y_1) = p(y_2)$ ; 
thus, $ y_2 = y_1$, so $(y_1, j_1) = (y_2, j_2)$.    
\end{enumerate}

\vspace{12pt}
\noindent
\underline{ $\tilde{p}|_{N(x, i)}$ is a surjection}.
Let $(u, k) \in N((p(x), i)$.  Again there are two cases to consider:
\begin{enumerate}
\item Suppose $u = p(x)$ and $i \sim k$.  Then $\tilde{p}(x, k) = (u, k)$, 
and by definition, $(x, k) \in N(x, i)$.
\item Suppose $u \sim p(x)$ and $i = k^{\phi(u p(x))}$.  Since $p|_{N(x)}$ is 
a surjection, there exists some $y \in N(x)$ such that $p(y) = u$.  Then 
$y \sim x$ and $i = k^{\phi(p(y)p(x))} = k^{{\psi}(yx)}$. 
Thus, $(y, k) \in N(x, i)$ and $\tilde{p}(y, k) = (u, k)$.  
\end{enumerate} 
Therefore, $\tilde{p}|_{N(x, i)}$ is a bijection.
This completes the proof.
\end{proof}

The vertex set of a bundle over $G$ is $V(G){\times}V(F)$. Then
$$L^2(G{\times}^{\phi}F) = L^2(G){\otimes}L^2(F)$$
where the tensor product takes place in the category of Hilbert spaces.
More precisely, it is the completion of the algebraic tensor product.
The following theorem (proved in \cite{kwak} for the finite case) 
provides a decomposition for the adjacency operator of any graph bundle.

\begin{thm}\label{adjacency}
Let $\phi$ be an $\text{Aut}(F)$-voltage assignment on a locally finite graph
$G$, with $F$ locally finite. Then 
$$A_{G \times^\phi F} = \bigoplus_{\gamma \in 
\mathrm{ Aut }(F)} A_{{{\overrightarrow{G}}}_{(\phi, \gamma)}} \otimes 
P_\gamma + I_G \otimes A_F.$$
\end{thm}

\begin{proof}
It is enough to prove the result for functions of the form $f{\otimes}g$,
where $f\in L^2(G)$ and $g\in L^2(F)$. Let
$(u, i) \in V(G{\times}^{\phi}F)$. Then
$$A_{G \times^\phi F}(f{\otimes}g)(u, i) = \sum_{(u, i) \sim (v, j)}f(v)g(j).$$
The right hand side is given by:
$$
\bigoplus_{\gamma \in 
\mathrm{ Aut }(F)} A_{{{\overrightarrow{G}}}_{(\phi, \gamma)}} \otimes 
P_\gamma(f{\otimes}g)(u, j) + I_G \otimes A_F(f{\otimes}g)(u, j)  = 
\bigoplus_{\gamma \in \mathrm{ Aut }(F),u\in
{{\overrightarrow{G}}}_{(\phi, \gamma)}} 
A_{{{\overrightarrow{G}}}_{(\phi, \gamma)}}(f|)(u)P_{\gamma}(g)(i)
+ f(u)A_F(g)(i)
$$
There are two possibilities for $(u, i) \sim (v, j)$:
\begin{enumerate}
\item $u \sim v$ and $i = j^{\phi(uv)}$. Then the right hand side  
becomes:
$$A_{{{\overrightarrow{G}}}_{(\phi, {\phi}(uv))}}(f|)(u)P_{\gamma}(g)(i)
= f(v)g(j).$$
\item $u = v$ and $i \sim j$. In the right hand side, only the last summand
is non-zero and it is equal to $f(u)g(j)$.
\end{enumerate}

Finally, it is clear that if neither $u \sim v$ nor $u = v$, then 
the sum on the right hand side is zero.
This completes the proof.
\end{proof}

By a {\it marked graph}, we mean a pair $(X, w)$ with $X$ a graph and $w$ a 
distinguished vertex.

\begin{defn}
On the space of marked graphs there is a metric $\text{dist}$ 
defined as follow:
$$\text{dist}\Big((X_1, w_1), (X_2, w_2)\Big) = 
\inf\left\{\frac{1}{n+1}; B_{X_1}(w_1, n) 
\text{is isometric to } B_{X_2}(w_2, n)\right\},$$
where $B_X(w, n)$ is the combinatorial ball of radius $n$ in $X$ centered 
on $v$.

For a sequence of marked graphs $\{(X_n, w_n)\}_{n\in\bbN}$, we say that 
$(X, w)$ is the limit graph if 
$$\lim_{n \to \infty} \text{dist}\left((X, w), (X_n, w_n)\right) = 0.$$
\end{defn}

For a finite graph $X$, let $c_r(X)$ denote the number of closed paths 
in $X$ of length $r$. Let
$$(X, w) = \lim_{n \to \infty} (X_n, w_n),$$ 
where $\{(X_n, w_n)\}_{n\in\bbN}$ is a covering sequence of 
$k$-regular marked graph. In \cite{grig}, the definition of the number 
$c_r$ is extended 
for the graph $X$ as follows:
$$\tilde{c}_r = \lim_{n \to \infty} \frac{c_r(X_n)}{|X_n|}.$$
In \cite{grig}, it was shown that the limit exists.
The zeta function $\zeta(X, w)$ of the marked graph $X$, with respect to 
the sequence $\{(X_n, w_n)\}_{n\in\bbN}$, is defined by 
$$\mathrm{ln }\zeta_{(X, w)}(z) = \lim_{n \to \infty} \frac{1}{|X_n|} 
\ln\zeta_{X_n}(z) = \sum_{r = 1}^\infty \frac{\tilde{c_r}z^r}{r}, \quad
|z| < \frac{1}{k - 1}.$$ 
The proof that the series has a non-trivial radius of convergence is given
in \cite{grig} and depends on results from \cite{serre}.

Let $X$ be a graph such that the degrees of vertices is bounded.
Let $\Gamma$ be a group of graph automorphisms of the graph $X$ that acts
on $X$ without inversions and satisfying the following properties:
\begin{enumerate}
\item For each $v\in V(X)$, the stabilizer ${\Gamma}_v = \{
{\gamma}\in{\Gamma}:\; {\gamma}v = v\}$ is finite.
\item If ${\F}_0 \subset V(X)$ is a complete set of orbit representatives
of the action of $\Gamma$ on $V(X)$, then
$$\text{vol}(X/{\Gamma}) = \sum_{v\in {\F}_0}\frac{1}{|{\Gamma}_v|} < 
{\infty}.$$
\end{enumerate}
In particular, if the action of $\Gamma$ on $V(X)$ is free, the second 
condition is equivalent to the condition that the orbit space $V(X)/{\Gamma}$
is finite. In this case, the Ihara zeta function is defined as 
$${\zeta}_{X,{\Gamma}}(z) = \prod_{C\in {\p}/{\Gamma}}\left(
1 - z^{{\ell}(C)}\right)^{|{\Gamma}_C|}.$$
where:
\begin{itemize}
\item $\p$ are equivalence classes of closed, primitive, tail-less edge-paths
without backtracking. Two such circuits are equivalent if they differ only by
a shift. ${\p}/{\Gamma}$ denotes the orbit space of $\p$ under the $\Gamma$
action.
\item For each class $C \in {\p}/{\Gamma}$, ${\ell}(C)$ denotes
the length of $C$ i.e., the number of edges in $C$.
\item ${\Gamma}_C$ denotes the isotropy group of $C$.
\end{itemize}
This formula generalizes the classical zeta function on finite graphs.

We will describe the analogue of Bass' formula for ${\zeta}_{X,{\Gamma}}(z)$
Let $L^2(X)$ be the Hilbert space of functions on $V(X)$. 
A unitary representation is given 
by:
$${\lambda}_0: {\Gamma} \to U(L^2(X)), \;\;
({\lambda}_0({\gamma})f)(v) = f({\gamma}^{-1}v), \;{\gamma}\in{\Gamma},\;
f\in L^2(X), \; v\in V(X).$$
Then the von Neumann algebra of all bounded operators on $L^2(X)$ that
commute with the $\Gamma$ action is defined as:
$${\N}_0(X, {\Gamma}) = \{{\lambda}_0({\gamma}): \; {\gamma}\in {\Gamma}\}'.$$
The algebra ${\N}_0(X, {\Gamma})$ inherits a trace given by:
$$\text{Tr}_{\Gamma}(A) = \sum_{v\in {\F}_0}\frac{1}{|{\Gamma}_v|}A(v, v),
\;\; A\in {\N}_0(X, {\Gamma}).$$
With this setting, the Bass formula for the Ihara zeta function has the form
(\cite{clair}, \cite{clair2}, \cite{guido1}, \cite{guido2}, \cite{guido3}):
$${\zeta}_{X,{\Gamma}}^{-1}(z) = 
(1 - z^2)^{-{\chi}^{(2)}(X)}\text{det}_{\Gamma}(
{\Delta}_{X,z}),$$
where 
\begin{itemize}
\item $\text{det}_{\Gamma} = \text{exp}{\circ}\text{Tr}_{\Gamma}{\circ}{\ln}$
is the determinant in the von Neumann algebra ${\N}_0(X, {\Gamma})$.
\item ${\Delta}_{X,z} = I - Az + Qz^2$, with $A$ the adjacency operator on
$X$, and $Q$ is the operator on $L^2(X)$ given by:
$$Q(f)(v) = (\text{deg}(v) - 1)f(v), \quad \text{for each}\;\; v\in V(X).$$
\end{itemize}

\begin{rem}

\null
\par\noindent
\begin{enumerate}
\item In \cite{clair}, \cite{clair2}, \cite{guido1}, \cite{guido2}, 
\cite{guido3},
it was shown that the function ${\zeta}_{X,{\Gamma}}$ is defined 
for sufficiently small $|u|$. More precisely, if $k$ is the maximum
degree of $X$, ${\zeta}_{X,{\Gamma}}(u)$ is a holomorphic function
for all $|u| < \frac{1}{d-1}$.
\item ${\chi}^{(2)}(X)$ is the Euler characteristic defined in
\cite{clair}. In most applications, it is equal to ${\chi}(X/{\Gamma})$,
the Euler characteristic of the orbit space.
\item Let $X$ be a $k$-regular graph and $q = k - 1$.
Using the determinant formula, the zeta function can be extended
to a holomorphic function in the open set (\cite{clair}, \cite{guido1}):
$${\Omega}_q = {\bbR}^2{\setminus}\left(
\left\{(x, y)\in{\bbR}^2:\; x^2 + y^2 = \frac{1}{q}\right\}\bigcup
\left\{(x, 0)\in{\bbR}^2:\; \frac{1}{q}\le |x| \le 1\right\}
\right).$$
\item In the above references there is an interpretation of the Bass
formula over the determinant on ${\N}_1(X, {\Gamma})$, the 
von Neumann algebra on the set of edges of $X$.
\end{enumerate}
\end{rem}

\vspace{12pt}\noindent
{\bf Notation}. There are three types of zeta functions used in this paper.
\begin{enumerate}
\item We write ${\zeta}_X(z)$ for the classical zeta function defined 
for a finite graph $X$.
\item We write ${\zeta}_{X,{\Gamma}}(z)$ for the equivariant zeta function
defined on an infinite graph $X$ equipped with an action of a group $\Gamma$
with finite co-volume.
\item We write ${\zeta}_{(X, w)}(u)$ for the zeta function that it is the 
limit of ${\zeta}_{X_n}(z)^{1/|V(X_n)|}$, where $\{(X_n, w_n)\}_{n\in\bbN}$
is a covering sequence of finite regular graphs converging to $(X, w)$.
\end{enumerate}

\vspace{12pt}
\begin{defn}
The sequence $\{X_n, w_n\}_{n\in\bbN}$ strongly converges to $(X, w)$ if:
\begin{enumerate}
\item $\{X_n, w_n\}_{n\in\bbN}$ is a covering sequence of marked 
$k$-regular graphs with
$$p_{m-1}: X_m \to X_{m-1}$$
the covering map.
\item $X$ is $k$-regular.
\item There are covering maps
$$\rho_n : X \to X_n$$
such that:
\begin{enumerate}
\item $\rho_1(u) = p_1p_2...p_{n-1}(\rho_n(u))$, for all $n$.
\item For each $n$, the isometry between $B_X(w, s_n)$ and 
$B_{X_n}(w_n, s_n)$ is given by the restriction of ${\rho}_n$.
\end{enumerate}
\end{enumerate}
\end{defn}

\begin{rem}
Cayley graphs of groups give sequences of graphs that strongly converge. Let $\Gamma$ be a group, $S$ a symmetric
set of generators and $\{K_n\}_{n\in\bbN}$ a sequence of normal
subgroups of finite index such that:
$$K_1 \supset K_2 \supset K_3 \dots , \quad 
\text{and} \quad \bigcap_{n\in\bbN}K_n = \{K\}.$$
Then the sequence of the marked Schreier graphs 
$\{({\s}({\Gamma}, K_n, S), 1)\}_{n\in\bbN}$
strongly converges to $({\s}({\Gamma}, K, S), 1)$.
\end{rem}

Let $\{X_n, w_n\}_{n\in\bbN}$
strongly converge to $(X, w)$
and ${\Gamma}_n = \text{Cov}(X, X_n)$. The next results gives a connection
between the different types of zeta functions.

\begin{thm}\label{thm-limit}
Assume that all the graphs in the sequence are $k$-regular finite
graphs converging to a $k$-regular graph
$(X, w)$. Then
$$\lim_{n\to \infty}{\zeta}_{X_n}(z)^{\frac{|V(X_1)|}{|V(X_n)|}} 
= {\zeta}_{X,{\Gamma}_1}(z) = {\zeta}_{(X,w)}(z)^{|V(X_1)|}.$$
\end{thm}

\begin{proof}
By \cite{clair2}, Theorem 2.1,
$$\lim_{n\to \infty}{\zeta}_{X_n}(z)^{\frac{1}{|N_n|}} 
= {\zeta}_{X,{\Gamma}_1}(z)$$
where $N_n = [{\Gamma}_n, {\Gamma}_1] = |V(X_n)|/|V(X_1)|$. The
result follows from the definition of ${\zeta}_{(X,w)}(z)$.
\end{proof}

The following is the main part of the proof of Theorem 2.1 in \cite{clair2}.

\begin{cor}\label{cor-determinant}
With the above notation,
$${\det}_{{\Gamma}_1}({\Delta}_{X,z}) = \lim_{n\to\infty}\left(
{\Delta}_{X_n,z}\right)^{\frac{|V(X_1)|}{|V(X_n)|}}$$
\end{cor}

Let $F$ be a locally finite graph and ${\phi}_1$ an $\text{Aut}(F)$-voltage 
assignment on $X_1$. Inductively, define an $\text{Aut}(F)$-voltage 
assignment on $X_n$ by:
$$\phi_n(uv) = \phi_{n-1}(p_{n-1}(u)p_{n-1}(v)).$$ 
 Also, define an $\text{Aut}(F)$-voltage 
assignment on $X$ by:
$$\phi(uv) = \phi_1(\rho_1(u)\rho_1(v)).$$
The details are presented in the following diagram.

$$
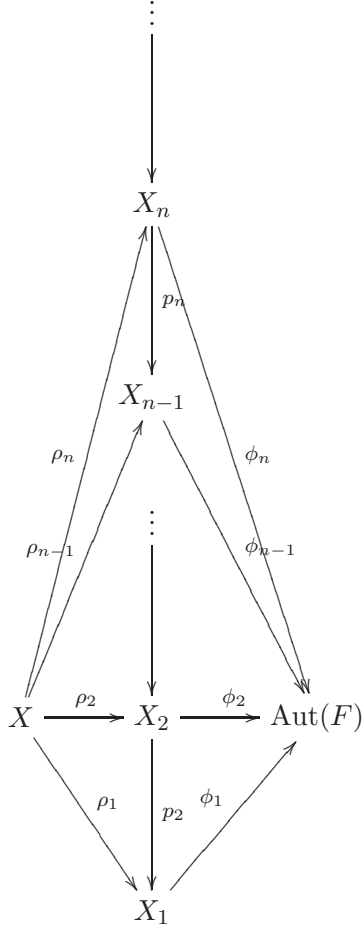
\xymatrix{
& \vdots \ar[dd]& \\
& & \\
& X_n \ar[dddddr]^{{\phi}_n} \ar[dd]^{p_n} & \\
& & \\
& X_{n-1} \ar[dddr]^{{\phi}_{n-1}} \\
& \vdots \ar[dd] & \\
& & \\
X \ar[uuuuur]^{{\rho}_n} \ar[uuur]^{{\rho}_{n-1}}
\ar[r]^{{\rho}_2} \ar[ddr]^{{\rho}_1} 
&X_2 \ar[r]^{{\phi}_2} \ar[dd]^{p_2}& \text{Aut}(F) \\
& & \\
& X_1 \ar[uur]^{{\phi}_1} &
}$$

Now, by Theorem \ref{cover}, we know that, for any finite, $d$-regular graph 
$F$, the sequence $\{X_n \times^{\phi_n} F\}_{n\in\bbN}$ is a $kd$-regular 
covering sequence; thus, by \cite{grig}, it converges.  We will show that 
in fact it converges to the graph $X \times^{\phi} F$.  
To do this will need the following:

\begin{lem}{\label{infinitecover}} 
Assume that $\{(X_n, w_n)\}_{n\in\bbN}$ strongly converges to
$(X, w)$. Then 
$$\widetilde{\rho_n}: B_{X \times^\phi F}((w, i), s_n) \to 
B_{X_n \times^{\phi_n} F}((w_n, i), s_n), \quad 
\widetilde{\rho_n}(u, i) = ((\rho_n(u), i),$$ 
is an isometry, for any $i \in V(F)$ and for all $n\in \bbN$.  
\end{lem}

\begin{proof}
Since $\rho_n$ is a bijection, it is clear that $\widetilde{\rho_n}$ is a 
bijection; thus, we must show that $\widetilde{\rho_n}$ preserves adjacency.  
To this end, assume $(u, i) \sim (v, j)$, for 
$(u, i), (v, j) \in B_{X \times^\phi F}(w, i)$.  Then there are two cases:
\begin{enumerate}
\item $u \sim v$ and $j = i^{\phi(uv)}$
\item $u = v$ and $i \sim j$.
\end{enumerate}
In the case where $u = v$, since clearly $\rho_n(u) = \rho_n(v)$, 
we must have $\widetilde{\rho_n}(u, i) \sim \widetilde{\rho_n}(v, j)$.    
In the case where $u \sim v$, we must have $\rho_n(u) \sim \rho_n(v)$.  
Thus we must show that 
$$j = i^{\phi_n(\rho_n(u)\rho_n(v))}.$$
Now, by the definition of $\phi$, 
$$\phi(uv) = \phi_1(\rho_1(u)\rho_1(v)),$$ 
and by assumption, 
$$\phi_1(\rho_1(u)\rho_1(v)) = 
\phi_1(p_{n-1} ... p_1(\rho_n(u)), p_{n-1} ... p_1(\rho_n(v))).$$ 
But by the definition of $\phi_n$, 
$$\phi_n(\rho_n(u)\rho_n(v)) = \phi_1(p_{n-1} ... p_1(\rho_n(u)), 
p_{n-1} ... p_1(\rho_n(v))).$$
This shows that $$\phi_n(\rho_n(u)\rho_n(v)) = \phi(uv),$$
and thus $$j=i^{\phi(uv)} \Longrightarrow j=i^{\phi_n(\rho_n(u)\rho_n(v))}.$$
This shows that $\hat{\rho_n}$ preserves adjacency, and thus is an isometry.  
This completes the proof.
\end{proof}

As a corollary, we have the following theorem.

\begin{thm}\label{limit}
For each $i\in F$, the covering sequence 
$\{(X_n \times^{\phi_n} F, (w_n, i))\}_{n\in\bbN}$ strongly
converges to $(X \times^{\phi} F, (w, i))$.
\end{thm}

\begin{proof}
Theorem \ref{cover} implies that
the covering conditions of the strong convergence are satisfied 
The rest of the proof follows from Lemma \ref{infinitecover} and \cite{grig}.
\end{proof} 

\section{Zeta Functions for Bundles and Coverings}

In this section we will use our previous result and \cite{clair}, 
\cite{clair2}, \cite{guido1}, \cite{guido2}, \cite{guido3}, 
to generalize the results of \cite{feng} to infinite graph bundles.

\begin{defn}
Let $X$ be a graph equipped with an $\text{Aut}(F)$-voltage assignment
$\phi$. 
\begin{enumerate}
\item An action of a group $\Gamma$ on $X$ without edge inversions is called
$(F, {\phi})$-compatible if  
$${\phi}({\gamma}(u){\gamma}(v)) = {\phi}(uv), \;\;\text{for all}\;\;
u, v\in V(X), \;\gamma\in \Gamma .$$
\item An action of a group $\Delta$ without inversions on $F$ is 
called $(X, {\phi})$-compatible if $\text{Im}({\phi}) \subset
C_{\text{Aut}(F)}({\Delta})$ i.e., the image of $\phi$ centralizes
$\Delta$. 
\item The pair of groups $({\Gamma}, {\Delta})$ as before is called
$\phi$-compatible if $\Gamma$ is $(F, {\phi})$-compatible and 
$\Delta$ is $(X, {\phi})$-compatible.
\end{enumerate}
\end{defn}

\begin{lem}\label{compatible}
With the above notation, if the pair $({\Gamma}, {\Delta})$ is 
$\phi$-compatible, then the product action:
$$({\gamma}, {\delta})(x, i) = ({\gamma}x, {\delta}i),\;\;
({\gamma}, {\delta})\in {\Gamma}{\times}{\Delta},\;\;
(x, i)\in V(X{\times}^{\phi}F),$$
is an action by graph automorphisms on $X{\times}^{\phi}F$.
Furthermore, if the action of $\Gamma$ and $\Delta$ are of finite
co-volume, so is the action of ${\Gamma}{\times}{\Delta}$ on
$X{\times}^{\phi}F$.
\end{lem}

\begin{proof}
The proof follows from the definitions.
\end{proof}

\begin{thm}{\label{indigo}}
Assume that $({\Gamma}, {\Delta})$ is a pair of $\phi$-compatible actions.
Also, assume that the actions are of finite co-volume.  Then 
$${\zeta}_{X \times^\phi F,{\Gamma}{\times}{\Delta}}(z)^{-1} = 
(1 - z^2)^{-\chi^{(2)}(X \times^\phi F)} 
\; {\text{det}}_{{\Gamma}{\times}{\Delta}} \left( I  - 
\sum_{\gamma \in \text{Aut}(F)}(A_{{\overrightarrow{X}}_{(\phi, \gamma)}} 
\otimes P_\gamma + I_X \otimes A_F)z + Qz^2\right),$$ 
 $\chi^{(2)}$ is the Euler characteristic and $Q$ is the diagonal
operator with $(x, i)$-entry $\text{deg}(x) + \text{deg}(i) - 1$.
Furthermore, the zeta function is holomorphic for $|z| < \frac{1}{k + d - 1}$.
If $X$ is $k$ regular and $F$ $d$-regular, then 
${\zeta}_{X \times^\phi F,{\Gamma}}(z)$ can be extended to a 
holomorphic function on ${\Omega}_{k+d-1}$.
\end{thm} 

\begin{proof}
From \cite{guido1}, \cite{guido2}, we have 
$${\zeta}_{X\times^\phi F,{\Gamma}{\times}{\Delta}}(z)^{-1} = 
(1 - z^2)^{-\chi^{(2)}(X \times^\phi F)} \; {\text{det}}_{{\Gamma}{\times}
{\Delta}} 
\left( I - zA_{X \times^{\phi} F} + z^2Q_{X \times^\phi F} \right).$$
The theorem now follows immediately from Theorem \ref{adjacency}.

\end{proof}

We will now use Theorem \ref{indigo} to provide a decomposition for the 
zeta function of any infinite cover. Let $p: Y \to X$ 
be a cover with $X$ finite and $Y$ locally finite. Let 
$\text{Cov}(p) = \Gamma$.
Now we define the function 
$$\phi: E({\overrightarrow{X}}) \to \Gamma.$$  
For this we write $X = \{x_1, \ldots x_n\}$.  
For each $i$, choose $v_i \in Y$ such that $p(v_i) = x_i.$  
Now, since 
$p: N(v_i) \to N(x_i)$ is a bijection, for each 
$x_j \in N(x_i)$ there exists a unique $u_j \in N(v_i)$ such that 
$p(u_j) = x_j$. So, since 
$p(v_j) = x_j = p(u_j)$,
there exists some $\gamma \in \Gamma$ such that 
$\gamma v_j = u_j$. Thus, define
$$\phi: E({\overrightarrow{X}}) \to 
\Gamma, \quad  \phi(x_ix_j) = \gamma.$$ 
We then have the following:

\begin{lem}{\label{folklore}}
Let $\phi$ be defined as above.  Then
\begin{enumerate}
\item The $\Gamma$ action on $Y$ is of finite co-volume and it 
is $\phi$-compatible.  
\item The map
$$\alpha: Y \to X 
\times^\phi \Gamma_, \quad \alpha(u) = \big(p(u), \beta \big)$$
is an isomorphism, where $\beta u = v_i, \quad p(u) = x_i = p(v_i)$.
\end{enumerate}
\end{lem}

\begin{proof}
The proof is folklore.  
\end{proof}

Now, in order to prove an analogue of \ref{adjacency} for 
$Y \simeq X \times^\phi \Gamma_1$, 
we need to define the following operator: 
for $\gamma_1, \gamma_2, \gamma \in \Gamma_1$, 
$$P_\gamma(\gamma_1, \gamma_2) = 
\begin{cases}
1,  \text{ if } \gamma_2 = \gamma_1\gamma \\
0,  \text{ otherwise} . 
\end{cases}
$$

\begin{lem}
Let $\phi$ be defined as above.  Then 
$$A_{X \times^\phi \Gamma} = 
\sum_{\gamma \in \Gamma} 
A_{{\overrightarrow{X}}_{(\phi, \gamma)}} \otimes P_\gamma.$$
\end{lem}

\begin{proof}
The proof is analogous to that of \ref{adjacency}. 
\end{proof}

The following theorem provides a decomposition for the zeta function 
of any infinite cover.

\begin{thm}
With the above notation,
$${\zeta}_{Y,{\Gamma}}(z)^{-1} = (1 - z^2)^{-\chi^{(2)}(X)} \; 
{\text{det}}_{\Gamma} \left(I - \left(\sum_{\gamma \in \Gamma} 
A_{{({\overrightarrow{X}}_{(\phi, \gamma)}})} \otimes 
P_\gamma\right)z + Qz^2 \right),$$ 
where $Q$ is the diagonal operator with $(x, {\gamma})$-entry
$\text{deg}(x) - 1$. The function is holomorphic for $|z| < \frac{1}{k-1}$.
If $Y$ is $k$-regular then ${\zeta}_{Y,{\Gamma}}(z)$
is holomorphic on ${\Omega}_{k-1}$.
\end{thm}

\begin{proof}
This follows immediately from \ref{indigo} and \ref{folklore}.
\end{proof}

Let $\{(X_n, w_n)_{n\in\bbN}$ be a sequence of finite $k$-regular marked graphs
that strongly converges to the $k$-regular marked graph $(X, w)$. Let
$F$ be a finite $d$-regular graph. With the notation as in Theorem 
\ref{limit}, we know that

We write $a_n = |V(X_n)|$ and $f = |F|$.

\begin{cor}\label{cor-limit}
With the above notation, for $|z| < \frac{1}{k + d - 1}$,
$$\begin{array}{ll}
&{\zeta}_{(X{\times}^{\phi}F, (w, i))}(z)^{-1}  = 
 \\[1ex]
 & =  \displaystyle{\lim_{n\to\infty}\left[
(1 - z^2)^{-{\chi}(X_n \times^{{\phi}_n} F)} 
\; {\text{det}}\left( I  - 
\sum_{\gamma \in \mathrm{Aut}(F)}(A_{{\overrightarrow{X_n}}_{({\phi}_n, 
\gamma)}} 
\otimes P_\gamma + I_{X_n} \otimes A_F)z + Q_nz^2\right)
\right]^{\frac{1}{fa_n}}}
\\[1ex]
 & =  \displaystyle{\left[
(1 - z^2)^{-\chi^{(2)}(X \times^\phi F)} 
\; {\text{det}}_\Gamma \left( I  - 
\sum_{\gamma \in \mathrm{Aut}(F)}(A_{{\overrightarrow{X}}_{(\phi, \gamma)}} 
\otimes P_\gamma + I_X \otimes A_F)z + Qz^2\right)\right]^{\frac{1}{fa_1}}}
\end{array}$$
\end{cor}

\begin{proof}
The first identity follows  because of Theorem \ref{limit}:
$$\begin{array}{ll}
&{\zeta}_{(X{\times}^{\phi}F, (w, i))}(z)^{-1}  = 
\displaystyle{
\lim_{n\to \infty}{\zeta}_{X_n{\times}^{{\phi}_n}F}(z)^{-\frac{1}{fa_n}}} 
 \\[1ex]
 & =  \displaystyle{\lim_{n\to\infty}\left[
(1 - z^2)^{-{\chi}(X_n \times^{{\phi}_n} F)} 
\; {\text{det}}\left( I  - 
\sum_{\gamma \in \text{Aut}(F)}(A_{{\overrightarrow{X_n}}_{({\phi}_n, 
\gamma)}} 
\otimes P_\gamma + I_{X_n} \otimes A_F)z + Q_nz^2\right)
\right]^{\frac{1}{fa_n}}}
\end{array}$$
The second identity follows from Theorem \ref{indigo} and Theorem 
\ref{thm-limit}.
\end{proof}

\section{Application}

Let $p: Y \to X$ be a cover with $\textbf{Cov}(p) = {\Gamma}$ and
$X$ a finite graph. let $F$ be a finite $d$-regular graph with
$n$ such that $\text{Aut}(F)$ contains a the dihedral group $D_{2n}$ of
order $2n$. Let $\phi$ be an $\text{Aut}(F)$-voltage assignment on $X$ 
whose image is contained into $D_{2n}$ and
$\psi$ the induced $\text{Aut}(F)$-voltage assignment on $Y$ (Theorem
\ref{cover}). By Theorem \ref{cover}, the induced map
$$\tilde{p}: Y{\times}^{\psi}F \to X{\times}^{\phi}F, \quad
\tilde{p}(x, i) = (p(x), i)$$
is a covering map. Also, $\text{Cov}(\tilde{p}) = {\Gamma}$. 

The following is the setup (for the finite case this is the same as in 
\cite{kwakk} and \cite{feng}): set 
$V(F) = \{1, 2, \ldots, n\}$ and $S_n$ the symmetric
group on $V(F)$.
Let $a = (1\;\; 2\; \ldots \; n-1 \;\; n)$ be an $n$-cycle and let 
$$b = \left\{\begin{array}{ll}
\displaystyle{(1 \;\; n)(2 \;\; n-1) \ldots \left(\frac{n-1}{2} \;\;
\frac{n+3}{2}\right)
\left(\frac{n+1}{2}\right)
}&\text{if $n$ is odd,} \\[1ex]
\displaystyle{
(1 \;\; n)(2 \;\; n-1) \ldots \left(\frac{n}{2}\;\; \frac{n+2}{2}\right)}& 
\text{if $n$ is even}
\end{array}\right.$$
be a permutation in $S_n$.  The permutations $a$ and $b$ generate the 
dihedral subgroup $D_n$ of $S_n$:  
$$D_n = {\langle} a, b | \; a^n = b^2 = 1, \; bab = a^{-1}{\rangle}.$$
Let $\mu = \mathrm{exp}(2\pi i / n)$ and ${\bf x}_k =
(1, {\mu}^k, {\mu}^{2k}, \dots , {\mu}^{(n-1)k})^T$ be the column
vector in ${\bbC}^n$. Then $1$, ${\mu}$, \dots  , ${\mu}^{n-1}$ are
the distinct eigenvalues of the permutation matrix $P(a)$ and ${\bf x}_k$
is the eigenvector corresponding to the eigenvalue ${\mu}^k$. Let $P(b)$
be the permutation matrix of $b$ and
$$M = \left\{
\begin{array}{lllllllll}
[{\bf x}_0 & P(b){\bf x}_1 & x_2 & P(b){\bf x}_2 & \dots & {\bf}x_{(n-1)/2} 
& P(b){\bf}x_{(n-1)/2}] && \text{if $n$ is odd} \\[1ex]
[{\bf x}_0 & P(b){\bf x}_1 & x_2 & P(b){\bf x}_2 & \dots & {\bf}x_{(n-2)/2} 
& P(b){\bf}x_{(n-2)/2} & {\bf x}_{n/2}] & \text{if $n$ is even}
\end{array}\right.
$$
In \cite{kwakk}, (also \cite{feng}) 
it was shown that $P(b){\bf x}_k$ is an eigenvector of $P(a)$ 
associated with the eigenvalue ${\mu}^{n-k}$. Thus $M$ is invertible. Also,
$P(a)$ and $A_F$ commute and thus they are simultaneously diagonalizable.
Also, ${\bf x}_k$ and $P(b){\bf x}_k$ ($1 \le k \le (n - 1)/2$ when $n$ is odd
and $1 \le k \le (n - 2)/2$ when $n$ is even)
are eigenvectors of $A_F$ with
the same eigenvalue of $P(b)$, denoted ${\lambda}_{(F,k)}$
Also, ${\bf x}_0$ is the eigenvector of $A_F$ corresponding to the
eigenvalue $d$ and, for $n$ even, ${\lambda}_{(F, n/2)}$ is the
eigenvalue associated to the eigenvector ${\bf x}_2$. Then as in \cite{kwakk},
using Theorem \ref{adjacency}, we get that 
$$\begin{array}{ll}
& (I_Y{\otimes}M)^{-1}A_{Y{\times}^{\psi}F}(I_Y{\otimes}M) = \\[2ex]
& = \displaystyle{
\left\{\begin{array}{ll}
(A_Y + dI_Y){\oplus}\left(\bigoplus_{i=1}^{(n-1)/2}(A_t + {\lambda}_{(F,t)}
(I_Y{\oplus}I_Y)
\right) & \text{if $n$ is odd}\\[5ex]
(A_Y + dI_Y){\oplus}\left(\bigoplus_{i=1}^{(n-2)/2}(A_t + {\lambda}_{(F,t)}
(I_Y{\oplus}I_Y) {\oplus}(B + {\lambda}_{(F,n/2)}I_Y
\right) & \text{if $n$ is even}
\end{array}\right.}
\end{array}$$
where
$$B = \sum_{k=0}^{n-1}\left((-1)^kA(\overrightarrow{Y}_{(\psi, a^k)}) +  
(-1)^{k+1}A(\overrightarrow{Y}_{(\psi, a^kb)})\right),$$
and 
$$A_t = \sum_{k = 0}^{n-1}
\left[\begin{array}{ll}
\displaystyle{{\mu}^{tk}A(\overrightarrow{Y}_{(\psi, a^k)})} &
\displaystyle{{\mu}^{tk}A(\overrightarrow{Y}_{(\psi, a^kb)})} \\[3ex]
\displaystyle{{\mu}^{(n-t)k}A(\overrightarrow{Y}_{(\psi, a^kb)})} &
\displaystyle{{\mu}^{(n-t)k}A(\overrightarrow{Y}_{(\psi, a^k)})}
\end{array}\right]
$$
Also, let $L_{Y} = (Q_{Y} + dI_{Y}) \otimes I_2$.  
Then the calculation in section 4 in \cite{feng} can be carried through
in our setting and we get the following: 

\begin{thm}\label{dihedral}
Let $p: Y \to X$ be as above.
Then 
$$\zeta_{Y \times^\psi F,{\Gamma}}(z)^{-1} = 
(1 - z^2)^{-{\chi}^{(2)}(Y{\times}^{\psi}F)}f_{Y, F}(z) 
\prod_{t=1}^{(n-1)/2} g_{Y, F, t}(z)$$
when $n$ is odd, and
$$\zeta_{Y \times^\psi F, {\Gamma}}(z)^{-1} = 
(1 - z^2)^{-{\chi}^{(2)}(Y{\times}^{\psi}F)}
f_{X, F}(z) h_{Y, F}(z) \prod_{t=1}^{(n-2)/2} g_{Y,F,t}(z)$$
when $n$ is even, where
\begin{enumerate}
\item $g_{Y, F, t}(z) = 
{\det}_{\Gamma}\left(I_Y{\oplus}I_Y - 
\left(A_t + \lambda_{(F, t)} (I_Y{\oplus}I_Y)\right)z + L_{Y}z^2\right)$
\item $h_{Y, F}(z) =  
{\det}_{\Gamma}\left(I_Y - \left(B + \lambda_{(F, n/2)} I_Y\right)z + 
\left(Q_Y + d I_Y\right)z^2\right)$
\item $f_{X, F}(z) = {\det}_{\Gamma}\left(I_Y - \left(A_Y + dI_Y\right)z
+ \left(Q_Y + dI_Y\right)z^2
\right)$
\end{enumerate}
\end{thm}

\begin{proof} This follows by simple calculation from \ref{indigo} with 
$\Delta = \{1\}$, and 
Theorem 9 of \cite{feng}.
\end{proof}

Let $\{(X_m, w_m)\}_{m\in\bbN}$ be a sequence of finite regular graphs
that strongly converges to $(X, w)$. Let $F$ be a finite $d$-regular
with $n$ vertices such that $\text{Aut}(F)$ contains $D_{2n}$. Let
$\phi$ be an $\text{Aut}(F)$-voltage assignment on $X_1$ whose
image is contained in $D_{2n}$. Let ${\phi}_n$ be the induced 
$\text{Aut}(F)$-voltage on $X_m$ and $\psi$ be the induced 
$\text{Aut}(F)$-voltage assignment on $X$. Set
$${\Gamma} = \text{Cov}(X \to X_1), \quad
{\Delta}_m = \text{Cov}(X_m \to X_1), \;\; m\in\bbN.$$

\begin{cor}\label{cor-dihedral}
Let $a_m = |V(X_m)|$. With the above notation,
\begin{enumerate}
\item If $n$ is odd:
$${\zeta}_{X{\times}^{\psi}F,{\Gamma}}(z)^{-1} = 
{\zeta}_{(X{\times}^{\psi}F, (w, i))}(z)^{-a_1} = 
\lim_{m\to\infty}\left(
(1 - z^2)^{-{\chi}(X_m{\times}^{{\phi}_m}F)}f_{X_m, F}(z) 
\prod_{t=1}^{(n-1)/2} g_{X_m, F, t}(z)\right)^{\frac{a_1}{a_m}}.$$
\item If $n$ is even
$${\zeta}_{X{\times}^{\psi}F,{\Gamma}}(z)^{-1} = 
{\zeta}_{(X{\times}^{\psi}F, (w, i))}(z)^{-a_1} = 
\lim_{m\to\infty}\left(
(1 - z^2)^{-{\chi}(X_m{\times}^{{\phi}_m}F)}f_{X_m, F}(z)h_{X_m,F}(z) 
\prod_{t=1}^{(n-1)/2} g_{X_m, F, t}(z)\right)^{\frac{a_1}{a_m}}.$$
\end{enumerate}
\end{cor}

\begin{proof}
It follows from Theorem \ref{thm-limit}, Theorem \ref{limit} and
Theorem \ref{dihedral}.
\end{proof}

In some cases, we can get a better description of the functions appearing in
the expression for the zeta function of the limit. Assume that all the
graphs $X_m$, $m \in \bbN$, and $X$ are $p$-regular. Following
\cite{grig}, for each $m\in\bbN$ set:
$${\sigma}_m = \sum\frac{{\delta}_{{\lambda}_i(X_m)}}{a_m}$$
where ${\lambda})i(X_m)$ are the eigenvalues of the Markov operator
$(1/k)A_{X_m}$ on $X_m$ and ${\delta}_x$ is the Dirac function. The 
sequence $\{{\sigma}_m\}_{m\in\bbN}$ weakly converges to the spectral
measure $\sigma$ associated to the Markov operator $(1/k)A_X$. Using
the calculation of Section 5 in \cite{grig} and Corollary
\ref{cor-determinant}, we get:
$$\begin{array}{lll}
{\ln}f_{X,F}(z) & = & {\ln}{\det}_{\Gamma}\left(I_Y - \left(A_Y + dI_Y\right)z
+ \left(Q_Y + dI_Y\right)z^2\right) \\[2ex]
& = & \displaystyle{
\lim_{m\to\infty}\frac{1}{a_m}{\ln}{\det}\left(I_{X_m} - 
\left(A_{X_m} + dI_{X_m}\right)z
+ \left(Q_{X_m} + dI_{X_m}\right)z^2\right)}\\[2ex]
& = & \displaystyle{
\lim_{m\to\infty}\int_{-1}^1{\ln}\left(
1 - (p{\lambda} + d)z + (p - 1 + d)z^2\right)d{\sigma}_m({\lambda})}\\[2ex]
& = &  \displaystyle{
\int_{-1}^1{\ln}\left(
1 - (p{\lambda} + d)z + (p - 1 + d)z^2\right)d{\sigma}({\lambda})}
\end{array}$$
where ${\sigma}$ is the spectral measure associated to $(1/p)A_X$.

Summarizing:

\begin{cor}\label{cor-f}
With the above notation,
$${\ln}f_{X,F}(z) = \int_{-1}^1{\ln}\left(
1 - (p{\lambda} + d)z + (p - 1 + d)z^2\right)d{\sigma}({\lambda}), \;\;
\text{for}\;\; |z| < \frac{1}{p + d - 1},$$
where $\sigma$ is the spectral measure associated to the regular
random walk on $X$.
\end{cor}

We give a specific example. The same method works for any group for which
the spectral measure is known.
Let $\Gamma$ be the Grigorchuk group (\cite{bg}, \cite{dela}, \cite{grig0},
\cite{grig}). Then $\Gamma$ can be represented as a subgroup of
automorphisms of the rooted binary tree. Let $P = \text{St}(1^{\infty})$ be
the stabilizer of the infinite sequence of 1's. Let $P_m$ be the stabilizer
of all the elements that start with $m$ 1's and it has finite index in 
$\Gamma$. Then 
$$P = \bigcup_{m = 1}^{\infty}P_m.$$
If $S = \{a, b, c, d\}$ be the standard set of generators of $\Gamma$, then
the Schreier graphs $\{\mathbb{S}_m = {\s}({\Gamma}, P_m, S)\}_{m\in\bbN}$ 
converge to
$\mathbb{S} = {\s}({\Gamma}, P, S)$. 
All the graphs have as the base point the identity 
coset. Then in \cite{grig}, Corollary 9.2, we have that
$$\begin{array}{ll}
{\ln}{\zeta}_{\mathbb{S}, P}(z) & = \displaystyle{
-3{\ln}(1 - z^2) - 
\int_{-1/2}^{0} \frac{\left(1-8xz + 7z^2\right)|1-4x|}{2\pi
\sqrt{x(2x-1)(2x+1)(1-x)}}dx}\\[2ex]
& -  \displaystyle{
\int_{1/2}^{1} \frac{\left(1- 8xz + 7z^2\right)|1-4x|}{2\pi
\sqrt{x(2x-1)(2x+1)(1-x)}}dx}
\end{array}$$
Let $F$ be a $d$-regular graph as in the beginning of the section and
$\phi$ an $\text{Aut}(F)$-voltage assignment whose image lies into $D_{2n}$.
Combining Theorem \ref{dihedral} and Corollary \ref{cor-f} we get.

\begin{thm}\label{limit-dihedral}
With the above notation, let $\psi$ be the $\text{Aut}(F)$-voltage
assignment on $\mathbb{S}$ and $R = P_1/ P$. Then 
$$\begin{array}{ll}
&{\ln}{\zeta}_{(\mathbb{S}{\times}^{\psi}F, (P, i))}(z)   \\[2ex]
& \displaystyle{
=\frac{1}{da_1}\left[
{-\chi^{(2)}(\mathbb{S} \times^\phi F){\ln}(1 - z^2) +
{\ln}{\det}}_R \left( I  - 
\sum_{\gamma \in 
\text{Aut}(F)}(A({\overrightarrow{\mathbb{S}}}_{(\psi, \gamma)}) 
\otimes P_\gamma + I_{\mathbb{S}} \otimes A_F)z + 
Qz^2\right)\right]}
\end{array}$$
Furthermore,
$$\zeta_{(\mathbb{S} \times^\psi F, (P, i)}(z)^{-1} =
(1 - z^2)^{-{\chi}^{(2)}({\mathbb{S}}{\times}^{\psi}F)}f_{{\mathbb{S}}, F}(z) 
\prod_{t=1}^{(n-1)/2} g_{{\mathbb{S}}, F, t}(z)$$
when $n$ is odd, and
$$\zeta_{(\mathbb{S} \times^\psi F, (P, i)}(z)^{-1} = 
(1 - z^2)^{-{\chi}^{(2)}({\mathbb{S}}{\times}^{\psi}F)}
f_{X, F}(z) h_{{\mathbb{S}}, F}(z) \prod_{t=1}^{(n-2)/2} g_{{\mathbb{S}},F,t}(z)$$
when $n$ is even, where $g$ and $h$ are as in Theorem \ref{dihedral} and
$${\ln}f_{{\mathbb{S}}, F}(z) = \int_{-\frac{1}{2}}^0
\frac{{\ln}(1 - (8x + d)z + (7 + d)z^2)|1 - 4x|}{2{\pi}
\sqrt{x(2x - 1)(2x + 1)(1 - x)}}dx +
\int_{\frac{1}{2}}^1
\frac{{\ln}(1 - (8x + d)z + (7 + d)z^2)|1 - 4x|}{2{\pi}
\sqrt{x(2x - 1)(2x + 1)(1 - x)}}dx$$
\end{thm}


\frenchspacing

\begin{thebibliography}{10}

\bibitem{bg}
L. Bartholdi, R. I. Grigorchuk, \emph{On the spectrum of Hecke type 
operators related to some fractal groups}, Proc. Steklov Inst. Math.,  
\textbf{231} (2000), 1--41. 

\bibitem{bass}
H. Bass, \emph{The Ihara--Selberg zeta function of a tree lattice},
Internat. J. Math. \textbf{3} (1992), 717--797.

\bibitem{clair}
B. Clair, S. Mokhtari-Sharghi, \emph{Zeta functions of discrete groups 
acting on trees}, J. Algebra \textbf{237} (2001), 591--430. 

\bibitem{clair2}
B. Clair, S. Mokhtari-Sharghi, \emph{Convergence of zeta functions of graphs},
Proc. Amer. Math. Soc. \textbf{130} (2002), 1881--1886.

\bibitem{dela}
P. de la Harpe, 
\emph{Topics in geometric group theory},
Chicago Lectures in Mathematics. University of Chicago Press, 
Chicago, IL, 2000.

\bibitem{feng}
R. Feng, J. Kwak, K. Kim, \emph{Zeta functions of graph bundles}, preprint.

\bibitem{grig0}
R. I. Grigorchuk, V. V.  Nekrashevich, V. I. Sushchanski\u\i,
\emph{Automata, dynamical systems, and groups}, 
 Proc. Steklov Inst. Math., \textbf{231} (2000), 128--203. 

\bibitem{grig}
R. Grigorchuk, A. {\v Z}uk, \emph{The Ihara zeta function of infinite graphs, 
the KNS spectral measure and integrable maps}, in:
``Random Walks and Geometry'', Proc. Workshop (Vienna 2001),
V. A. Kaimanovich {\it et al.}, eds., de Gruyter, Berlin 2004, 141--180.

\bibitem{grig2}
R. Grigorchuk, A. {\v Z}uk, \emph{
On the asymptotic spectrum of random walks on infinite families of graphs},  
Random walks and discrete potential theory (Cortona, 1997),  
188--204, Sympos. Math., XXXIX, Cambridge Univ. Press, Cambridge, 1999. 

\bibitem{guido1}
D. Guido, T. Isola, M. L. Lapidus, \emph{Ihara zeta functions for 
periodic simple graphs}, arXiv:math.OA/0605753, May 2006.

\bibitem{guido2}
D. Guido, T. Isola, M. L. Lapidus, \emph{Ihara zeta functions for 
periodic simple graphs and its approximation in the amenable case}, 
arXiv:math.OA/0608229, August 2006.

\bibitem{guido3}
D. Guido, T. Isola, M. L. Lapidus, \emph{A trace on fractal graphs and the
Ihara zeta function}, 
arXiv:math.OA/0608060, August 2006.

\bibitem{ihara}
Y. Ihara, \emph{On discrete subgroups of the two by two projective linear 
group over p-adic fields}, J. Math. Soc. Japan \textbf{18} (1966), 
219--235.

\bibitem{kwakk}
J. H. Kwak, Y. S. Kwon, \emph{Characteristic polynomials of graph bundles 
having voltages in a dihedral group}, Linear Algebra Appl.  \textbf{336}  
(2001), 99--118. 

\bibitem{kwak}
J. H. Kwak, J.  Lee, \emph{Characteristic polynomials of some graph bundles 
I{\!}I},
Linear and Multilinear Algebra  \textbf{32} (1992), 61--73. 

\bibitem{alex}
A. Lubotzky, \emph{Discrete Groups, Expander Graphs and Invariant Measures}

\bibitem{serre}
J.-P. Serre, \emph{R{\'e}partitions asymptotique des valeurs propres de
l' op{\'e}rateur de Hecke $T_p$}, J. Amer. Math. Soc. \textbf{10} (1997),
75--102.

\bibitem{stark1}
H. M. Stark, A. A. Terras, \emph{Zeta functions of finite graphs and
coverings}, Adv. Math. \textbf{121} (1996), 126--165.

\bibitem{stark2}
H. M. Stark, A. A. Terras, \emph{Zeta functions of finite graphs and
coverings I{\!}I}, Adv. Math. \textbf{154} (2000), 132--195.

\bibitem{stark3}
H. M. Stark, A. A. Terras, \emph{Zeta functions of finite graphs and
coverings I{\!}I{\!}I}, Adv. Math.

\end{thebibliography}
\end{document}